\documentclass{article}
\usepackage{arxiv}
\usepackage[utf8]{inputenc} 
\usepackage[T1]{fontenc}    
\usepackage{hyperref}       
\hypersetup{hidelinks}
\usepackage{url}            
\usepackage{booktabs}       
\usepackage{amsfonts}       
\usepackage{nicefrac}       
\usepackage{microtype}      
\usepackage{lipsum}
\usepackage{graphicx}
\usepackage{algorithm}
\usepackage{algpseudocode}
\algtext*{EndIf} 
\algtext*{EndFor}
\usepackage{tikz}
\usepackage{amsmath,amssymb,amsthm}
\usetikzlibrary{trees}
\graphicspath{ {./images/} }

\usepackage[
  backend=biber,
  style=numeric,
  maxbibnames=99
]{biblatex}

\newtheorem{definition}{Definition}[section]

\usepackage{graphicx}  
\title{\resizebox{\textwidth}{!}{Computing Tools for Translation-Invariant Total Orders}}

\author{
Yifan Jing \\
 College of LSA\\
 University of Michigan \\
 Ann Arbor, MI 48104 \\
 \texttt{ivaning@umich.edu} \\
\And
Emma Yicheng Lou \\
 College of LSA\\
 University of Michigan\\
 Ann Arbor, MI 48104 \\
 \texttt{emmaycl@umich.edu} \\
 \And
Shaotian Sun \\
 College of LSA\\
 University of Michigan\\
 Ann Arbor, MI 48104 \\
 \texttt{shaotian@umich.edu} \\
 \And
Yuxuan Wu \\
 College of LSA\\
 University of Michigan\\
 Ann Arbor, MI 48104 \\
 \texttt{owenwho@umich.edu} \\
}

\begin{document}
\maketitle

\begin{abstract}
We introduce \textit{TITO\_Explore}, a software package for representing and computing with Translation-Invariant Total Orders (TITOs). We define a canonical window notation for TITOs and design and implement algorithms for several computational tasks involving them. The package normalizes the window notation of a given TITO into its canonical form, computes its inversion set, compares the weak order between two TITOs, and computes the join of two specified TITOs. Our weak order comparison algorithm operates by partitioning the inversion sets into disjoint subsets, thereby breaking down the comparison problem into evaluations of paired subsets. The join algorithm uses an edge-weighted directed graph to represent inversions and converts the problem of finding the join into a weighted path problem in the graph.
\end{abstract}


\section{Introduction}
\newcommand{\TTot}{\mathrm{TTot}}
A central theme in mathematics is the study of order structures that interact naturally with algebraic operations. In this paper, we consider \emph{translation-invariant total orders} (TITOs) on the integers, introduced by Barkley and Speyer \cite{Barkley_2024,barkley2025extendedweakorderaffine,barkley2025affinetamarilattice}.

\begin{definition}
A \emph{translation-invariant total order} (TITO) with period n, where n is a fixed positive integer, is a total order \(\prec\) on \(\mathbb{Z}\) such that for all \(a,b\in\mathbb{Z}\),
\[
a \prec b \quad \Longleftrightarrow \quad a+n \prec b+n.
\]
The set of TITOs is denoted by \(\TTot_n\).
\end{definition}

TITOs are closely related to the combinatorics of Coxeter groups and root systems, especially to the affine symmetric group, the affine Coxeter group of type \(\widetilde A_{n-1}\) \cite{Barkley_2024}. In this setting, TITOs encode combinatorial information connected with weak order, torsion classes, and lattice structures arising in algebraic combinatorics and representation theory. More concretely, the weak order on TITOs forms a
complete lattice, and this lattice is indeed a sublattice of the weak order on all total
orders of \(\mathbb Z\) \cite{barkley2024affineextendedweakorder}.

The role of inversion sets is central in this picture. In the classical theory of Coxeter groups, weak order can be described by the inclusion of inversion sets \cite{bjorner2005combinatorics}.
For infinite Coxeter groups, however, inversion sets and joins become more subtle. Hohlweg and Labbé study inversion sets, infinite reduced words, and weak order in Coxeter groups, emphasizing the need to understand weak order relations through appropriate inversion set data \cite{hohlweg2016inversion}.
TITOs provide a concrete affine setting in which these ideas can be explored computationally.

Because their inversion sets may be infinite, TITOs are difficult to handle algorithmically. As a result, basic tasks such as comparing two TITOs in weak order and computing the join of two TITOs require efficient finite representations and algorithms.

To address these problems, we develop finite representations for inversion sets of TITOs and implement algorithms for normalization, weak order comparison, inversion set computation, and join computation. These algorithms are collected in the Python package \textit{TITO\_Explore}. The goal of this paper is to explain the mathematical ideas behind these algorithms, describe their implementation, and illustrate their behavior through examples. The source code for the package is available at 
\nolinkurl{https://github.com/emmaycl/TITO_Explore_package}.

\section{Preliminaries}

We begin by introducing the notation needed to describe TITOs and the weak order between them.

A TITO may be represented in either one-line notation or window notation. To define these notations, we first introduce the notion of a block. A subset \(I\subseteq\mathbb{Z}\) is said to be \emph{order-convex} with respect to \(\prec\) if for any \(a,c\in I\) and any \(b\in\mathbb{Z}\) satisfying
\[
a \prec b \prec c,
\]
we also have \(b\in I\).

\begin{definition}
Let \(\prec\) be a TITO. A \emph{block} of \(\prec\) is an order-convex subset \(I\subseteq\mathbb{Z}\) satisfying the following conditions:
\begin{enumerate}
    \item The restriction of \(\prec\) to \(I\) has neither a minimal element nor a maximal element.
    \item For any \(a,c\in I\), the interval
    \[
    \{\, b\in I \mid a \prec b \prec c \,\}
    \]
    is finite.
\end{enumerate}
The \emph{size} of a block \(I\) is the number of residue classes modulo \(n\) represented in \(I\).

A block \(I\) is called \emph{waxing} if \(x \prec x+n\) for all \(x\in I\), and \emph{waning} if \(x+n \prec x\) for all \(x\in I\).
\end{definition}

The window notation is a finite representation for a TITO, introduced by Barkley \cite{barkley2025extendedweakorderaffine}: Each $k$-window represents a block, where we take $k$ consecutive numbers from the TITO to form the block. k is the size of the block. If the block is waning, we add an underline to the window. By concatenating the window for each block in the TITO, we get the window notation for the TITO. 

For example, when \(n=2\), the window notation \([0][\underline{1}]\) corresponds to the order
\[
\cdots \prec -2 \prec 0 \prec 2 \prec \cdots \prec 3 \prec 1 \prec -1 \prec \cdots .
\]
Note that in this TITO, the first block contains only numbers in the residue class 0, so the first window is a 1-window, as there is one residue class involved. The same argument applies to the second window. 

To define \emph{weak order} on TITOs, consider the equivalence relation on pairs \((a,b)\) given by
\[
(a,b)\sim(a+n,b+n).
\]
We denote the equivalence class of \((a,b)\) again by \((a,b)\). A \emph{reflection index} is such an equivalence class with \(a<b\).

\begin{definition}
An \emph{inversion} of a TITO \(\prec\) is a reflection index \((a,b)\) such that \(a < b\) and \(b \prec a\). Let \(N(\prec)\) denote the set of inversions of \(\prec\), \(N(\prec)\) is called the inversion set of the TITO. The \emph{weak order} on \(\TTot_n\) is defined by
\[
\prec_1 \le \prec_2
\quad \Longleftrightarrow \quad
N(\prec_1)\subseteq N(\prec_2).
\]
\end{definition}

We now illustrate this definition in the case \(n=2\). Let
\[
\prec_1=[0,1], \qquad \prec_2=\underline{[0,-1]}, \qquad \prec_3=[0][1].
\]
Then their inversion sets are
\[
N(\prec_1)=\emptyset,
\]
\[
N(\prec_2)=\{(0,2),(0,4),\ldots\}\cup \{(1,3),(1,5),\ldots\}
\cup \{(0,1),(0,3),\ldots\}\cup \{(1,2),(1,4),\ldots\},
\]
and
\[
N(\prec_3)=\{(1,2),(1,4),\ldots\}.
\]
Hence
\[
\prec_1 \le \prec_3 \le \prec_2,
\]
since
\[
N(\prec_1)\subseteq N(\prec_3)\subseteq N(\prec_2).
\]

Here, we also present an extended notation for a reflection index. 
\begin{definition}
Given a reflection index $(a,b)$ for a TITO with period $n$, the star notation for the reflection index is defined as
\[
(a,b)^* := \{(a,b+kn) \mid k \in \mathbb{Z}_{\geq0}\}.
\]
\end{definition}

\begin{definition}
Let $\prec$ be a TITO with period $n$. The \emph{reflection table} of $\prec$ is an $n\times n$ table $T$, whose rows and columns are indexed by residue classes modulo $n$. The entry $T[a,b]$ records the reflection indices in $N(\prec)$ whose first coordinate is congruent to $a \pmod n$ and whose second
coordinate is congruent to $b \pmod n$.

More precisely, an entry $T[a,b]$ consists of a finite set of increments $d$ satisfying
\[
    d>0,
    \qquad
    a+d \equiv b \pmod n.
\]
Each such increment represents the reflection index
\[
    (a,a+d).
\]
If the entry is marked with a star, then it represents the corresponding
infinite family
\[
    (a,a+d)^*.
\]
Thus the reflection table gives a finite representation of the inversion set
by grouping reflection indices according to their residue classes modulo $n$.
\end{definition}

\begin{definition}
Let $\widetilde T_n$ denote the set of reflection indices for TITOs of period $n$.
For a subset $X \subseteq \widetilde T_n$, the \emph{closure} of $X$, denoted
$\overline X$, is the smallest subset of $\widetilde T_n$ containing $X$ with the
following property: whenever there exist representatives $a<b<c$ such that
$(a,b)$ and $(b,c)$ belong to $\overline X$, the reflection index $(a,c)$ also
belongs to $\overline X$.
\end{definition}

\begin{definition}
The \emph{join} of two TITOs $\prec_1$ and $\prec_2$, denoted
$\prec_1 \vee \prec_2$, is the smallest TITO in weak order whose inversion set
contains both $N(\prec_1)$ and $N(\prec_2)$. Equivalently, its inversion set is
\[
N(\prec_1\vee\prec_2)
=
\overline{N(\prec_1)\cup N(\prec_2)}.
\]
\end{definition}
Computationally, we obtain $N(\prec_1\vee\prec_2)$ by starting with $N(\prec_1)\cup N(\prec_2)$ and adding all reflection indices forced by transitivity. In particular, if $(a,b)$ and $(b,c)$ are present, then $(a,c)$
must also be present.

\section{Implementation Details}
\subsection{Normalization}
The normalization procedure defines a standard window notation for a TITO, which satisfies the following condition: Fix $n\in \mathbb{Z}_{>0}$. In each block of the window notation, the first number is an integer from $0$ to $n-1$ and belongs to the smallest residue class among all the residue classes in this block. \\
The TITO normalization procedure normalizes each block using Algorithm~\ref{alg:normalization} and concatenates the normalized blocks. 
\begin{algorithm}
\caption{Window Notation Normalization (with Waxing/Waning)}
\label{alg:normalization}
\begin{algorithmic}[1]
\Require A block $W = [x_1,\dots,x_m]$ with distinct residues modulo $n$, $m \le n$, and status $s \in \{0, 1\}$ (0 for waxing, 1 for waning).

\State Compute residues $r_i = x_i \bmod n$ for all $i = 1,\dots,m$
\State Let $i^* = \arg\min r_i$

\If{$s = 0$} \Comment{Waxing condition}
    \State Reorder the block:
    \[
    W' = [x_{i^*}, x_{i^*+1}, \dots, x_m, x_1+n, \dots, x_{i^*-1}+n]
    \]
\ElsIf{$s = 1$} \Comment{Waning condition}
    \State Reorder the block:
    \[
    W' = [x_{i^*}, x_{i^*+1}, \dots, x_m, x_1-n, \dots, x_{i^*-1}-n]
    \]
\EndIf

\State \Return $W'$
\end{algorithmic}
\end{algorithm}

\paragraph{Example.}
The following code demonstrates the normalization procedure for a TITO with
period $n=3$:

\begin{verbatim}
from TITO_Explore.types import TranslationInvariantTotalOrder
from TITO_Explore.normalization import process_tito_blocks

tito = TranslationInvariantTotalOrder(
    n=3,
    num_blocks=2,
    vectors=[[4, 3], [5]],
    waxing_waning=[0, 0],
)

normalized_tito = process_tito_blocks(tito)
print(normalized_tito)
\end{verbatim}

The output is:
\begin{verbatim}
[[0, 4], [2]]
\end{verbatim}

Thus the window notation
\[
    [[4,3],[5]]
\]
is normalized to
\[
    [[0,4],[2]].
\]

\subsection{Weak Order Comparison}

Given two normalized TITOs $\prec_1$ and $\prec_2$ of period $n$, we compare them by comparing their inversion sets. The goal is to determine whether
\[
N(\prec_1)=N(\prec_2), \qquad
N(\prec_1)\subsetneq N(\prec_2), \qquad
N(\prec_2)\subsetneq N(\prec_1),
\]
or whether the two inversion sets are incomparable.

The main difficulty is that inversion sets are infinite. To overcome this, we decompose each inversion set according to pairs of residue classes modulo $n$. Each residue class pair contributes a local comparison result, and the global weak-order comparison is then obtained by combining these local results. Since there are $\frac{n(n+1)}{2}$ residue class pairs in total, this yields an $O(n^2)$ comparison algorithm.

\begin{algorithm}
\caption{Weak-order comparison of two normalized TITOs}
\label{alg:weakordercompare}
\begin{algorithmic}[1]
\Require Two normalized TITOs $\prec_1,\prec_2$ of period $n$
\Ensure One of $=$, $<$, $>$, or $\times$

\If{$\prec_1=\prec_2$}
    \State \Return $=$
\EndIf

\State $S \gets \emptyset$

\For{$i=0$ to $n-1$}
    \State Append $\Call{CompareSameResidue}{i,\prec_1,\prec_2}$ to $S$
\EndFor

\For{$i=0$ to $n-1$}
    \For{$j=i+1$ to $n-1$}
        \State $C_1 \gets \Call{ClassifyDifferentResiduePair}{i,j,\prec_1}$
        \State $C_2 \gets \Call{ClassifyDifferentResiduePair}{i,j,\prec_2}$
        \State Append $\Call{CompareDifferentResidue}{C_1,C_2}$ to $S$
    \EndFor
\EndFor

\If{$\times \in S$}
    \State \Return $\times$
\ElsIf{$\subsetneq \in S$ and $\supsetneq \in S$}
    \State \Return $\times$
\ElsIf{$\subsetneq \in S$}
    \State \Return $<$
\ElsIf{$\supsetneq \in S$}
    \State \Return $>$
\Else
    \State \Return $=$
\EndIf
\end{algorithmic}
\end{algorithm}

We now describe the local comparison rules used in Algorithm~\ref{alg:weakordercompare}.

\paragraph{Step 1: Decomposition by residue class pairs.}
The inversion set $N(\prec)$ decomposes into disjoint subsets indexed by residue class pairs modulo $n$. There are $n$ same-residue-class pairs and $\frac{n(n-1)}{2}$ different-residue-class pairs, for a total of
\[
\frac{n(n+1)}{2}
\]
components. Since reflection indices arising from different residue class pairs are distinct, the weak-order comparison can be carried out independently on each component.

\paragraph{Step 2: Same residue class pairs.}
For a fixed residue class $i \pmod n$, the corresponding local inversion subset is either empty or full as $(i,i+n)^*$. Hence the local comparison is immediate. If both local subsets are empty, or both are full, then the local result is equality. If the local subset for $\prec_1$ is empty while that for $\prec_2$ is full, then the local result is $\subsetneq$. If the local subset for $\prec_1$ is full while that for $\prec_2$ is empty, then the local result is $\supsetneq$.

\paragraph{Step 3: Different residue class pairs.}
Fix residue classes $i,j\in\{0,1,\dots,n-1\}$ with $i<j$. For the local comparison on the residue class pair $(i,j)$, we first classify the configuration of the two residue classes inside the normalized TITO.

\subparagraph{Step 3.1: Different-block cases.}
If the residue classes $i$ and $j$ lie in different blocks, then exactly two cases can occur, namely Case 1, written $[i][j]$, and Case 2, written $[j][i]$.

\subparagraph{Step 3.2: Same-block cases.}
If the residue classes $i$ and $j$ lie in the same block, we extract the corresponding 2-block and normalize it. The normalized 2-block has the form
\[
[i,t] \qquad\text{or}\qquad \underline{[i,t]},
\]
where $t\equiv j \pmod n$.

Among these same-block configurations, there are two special cases. In the usual-order case, the local inversion subset is empty. In the reversing-order case, the local inversion subset is full. Excluding these two special cases, there are four general same-block cases: Case 3, given by $[i,t]$ with $i<t$; Case 4, given by $[i,t]$ with $i>t$; Case 5, given by $\underline{[i,t]}$ with $i<t$; and Case 6, given by $\underline{[i,t]}$ with $i>t$.

\subparagraph{Step 3.3: Comparison rule for different residue class pairs.}
Let $C_1$ and $C_2$ be the classifications obtained from $\prec_1$ and $\prec_2$ for the same residue class pair $(i,j)$. If one local configuration is the usual order and the other is not, then the usual-order inversion subset is strictly smaller. If one local configuration is the reversing order and the other is not, then the reversing-order inversion subset is strictly larger. If both local configurations are in usual order, then the local result is equality, and the same conclusion holds when both local configurations are reversing order.

For the six non-special cases, define
\[
\text{Group A}=\{\text{Case 2},\text{ Case 3},\text{ Case 6}\},
\qquad
\text{Group B}=\{\text{Case 1},\text{ Case 4},\text{ Case 5}\}.
\]
If $C_1$ and $C_2$ lie in different groups, then the local result is incomparable, denoted by $\times$. Within Group A, the inclusion order is
\[
\text{Case 3} \subsetneq \text{Case 2} \subsetneq \text{Case 6},
\]
whereas within Group B, the inclusion order is
\[
\text{Case 4} \subsetneq \text{Case 1} \subsetneq \text{Case 5}.
\]
If $C_1$ and $C_2$ fall into the same case, then the comparison is refined further. In Cases 1 and 2, the local inversion subsets are equal. In Cases 3 and 6, the larger value of $t$ gives the larger inversion subset. In Cases 4 and 5, the larger value of $t$ gives the smaller inversion subset.

\paragraph{Step 4: Combination rule.}
After all local comparisons are completed, each residue class pair contributes one of the four outcomes
\[
\subsetneq,\qquad \supsetneq,\qquad =,\qquad \times.
\]
These local outcomes are then combined globally. If any local result is $\times$, then the final result is incomparable. Likewise, if both $\subsetneq$ and $\supsetneq$ occur among the local results, then the final result is also incomparable. If every non-equality local result is $\subsetneq$, then $N(\prec_1)\subsetneq N(\prec_2)$. If every non-equality local result is $\supsetneq$, then $N(\prec_2)\subsetneq N(\prec_1)$. Finally, if every local result is $=$, then $N(\prec_1)=N(\prec_2)$.

\paragraph{Example.}
The following code compares two TITOs in weak order.

\begin{verbatim}
from TITO_Explore.types import TranslationInvariantTotalOrder
from TITO_Explore.comparison import compare_titos_overall, relation_symbol
tito1 = TranslationInvariantTotalOrder(
    n=2,
    num_blocks=1,
    vectors=[[0, 1]],
    waxing_waning=[0],
)
tito2 = TranslationInvariantTotalOrder(
    n=2,
    num_blocks=2,
    vectors=[[0], [1]],
    waxing_waning=[0, 0],
)
print(compare_titos_overall(tito1, tito2))
\end{verbatim}

The output is:
\begin{verbatim}
tito1 < tito2
\end{verbatim}

Thus the first TITO is strictly smaller than the second TITO in weak order.

\subsection{Inversion Set Computation} \label{findinversionset}
This algorithm follows the comparison algorithm from the previous section. 
The complete inversion set of a TITO can be partitioned into disjoint subsets by analyzing each pair of residue classes. For each pair of residue classes, we extract the sub-window notations and apply the normalization algorithm. Then, for all the seven cases below, the inversions can be finitely represented using the reflection index notation along with the star notation. 

Let $i$ and $j$ be two residue classes. We abuse notation and identify $i$ and $j$ with their unique representatives between $0$ and $n-1$. 

\begin{itemize}
    \item \textbf{Case 1: }

    Let $i<j$ be two different residue classes. Suppose $i$ and $j$ are in
    different blocks, and the block containing $i$ appears before the block
    containing $j$. The normalized window notation is $[i][j]$, with $i<j$. Then the
    inversion contributed by this pair is
    \[
        (j, i+n)^*.
    \]
    Note that here, since each block only has one element, whether the block is waxing or waning doesn't affect the reflection indices brought by this pair of residue classes. Same applies to the second case. 
    
    \item \textbf{Case 2: }

    Let $i<j$ be two different residue classes. Suppose $i$ and $j$ are in
    different blocks, and the block containing $i$ appears after the block
    containing $j$. The normalized window notation is $[j][i]$, with $i<j$. Then the
    inversion contributed by this pair is
    \[
        (i,j)^*.
    \] 

    \item \textbf{Case 3: }

    Let $i<j$ be two different residue classes. Suppose $i$ and $j$ are in the
    same waxing block. Let $t$ be the representative of the residue class $j$
    appearing in the normalized 2-block, so that $t\equiv j \pmod n$.
    Suppose the normalized 2-block is
    \[
        [i,t],
        \qquad i<t.
    \]
    If $t=j$, 
    then this pair contributes no inversions. Otherwise, the inversions
    contributed by this pair are
    \[
        (i,t-kn),\ \ldots,\ (i,t-2n),\ (i,t-n),
    \]
    where $k$ is the largest positive integer such that
    \[
        t-kn>i.
    \] Here, $t-kn$ is exactly $j$. Here we abuse the notation so that $i$ stands for both the residue class $i$ and the first number in the local window notation. 

   \item \textbf{Case 4: }

    Let $i<j$ be two different residue classes. Suppose $i$ and $j$ are in the
    same waxing block. Let $t$ be the representative of the residue class $j$
    appearing in the normalized 2-block, so that $t\equiv j \pmod n$.
    Suppose the normalized 2-block is
    \[
        [i,t],
        \qquad t<i.
    \]
    Then the inversions contributed by this pair are
    \[
        (t,i),\ (t+n,i),\ \ldots,\ (t+\ell n,i),
    \]
    where $\ell$ is the largest nonnegative integer such that
    \[
        t+\ell n<i.
    \]
   
    \item \label{case5} \textbf{Case 5: }

    Let $i<j$ be two different residue classes. Suppose $i$ and $j$ are in the same waning block. Let $t$ be the representative of the residue class $j$ appearing in the normalized 2-block, so that $t\equiv j \pmod n$.
    Suppose the normalized 2-block is
    \[
        \underline{[i,t]},
        \qquad i<t.
    \]
    The inversions arising from this case include real reflection indices, namely inversions between numbers of two different residue classes, and imaginary reflection indices, namely inversions between numbers of the same residue class. Here, we only consider the real reflection indices; the imaginary reflection indices are added to the reflection table as described in the last case.
    
    Then the real inversions contributed by this pair are
    \[
        (i-n,t)^*,
        (t,i+\ell n)^*,
    \]
    where $\ell$ is the smallest positive integer, if it exists, such that $i+\ell n>t$.

    Equivalently, it is 
    \[
        (i-n,t)^*,
        (t, t+n-j+i)^*.
    \]

    It is also worth noting that the inversions brought by Case 5 are the complement of those brought by Case 4. The same representation for the inversions brought by Case 5 can be easily derived using this complement relation as well.

    \item \textbf{Case 6: }

    Let $i<j$ be two different residue classes. Suppose $i$ and $j$ are in the
    same waning block. Let $t$ be the representative of the residue class $j$
    appearing in the normalized 2-block, so that $t\equiv j \pmod n$.
    Suppose the normalized 2-block is
    \[
        \underline{[i,t]},
        \qquad t<i.
    \]
    The inversions arising from this case include real reflection indices and
    imaginary reflection indices. Here, we only consider the real reflection
    indices.
    
    Then the real inversions contributed by this pair are
    \[
        (t,i)^*,
        \qquad
        (i,t+\ell n)^*,
    \]
    where $\ell$ is the smallest positive integer such that
    \[
        t+\ell n>i.
    \]

    Equivalently, it's 
    \[
        (t,i)^*,
        \qquad
        (i,j)^*.
    \]
    Similarly as in Case 5, note that the inversions brought by Case 6 are the complement of those brought by Case 3.

    \item \textbf{Imaginary Reflection Indices: }
    After all the real reflection indices in the inversion set are determined according to the first 6 cases, the algorithm checks all residue classes $0$ to $n-1$ and adds imaginary reflection indices brought by waning blocks. 
    
    Let $i$ be a residue class that appears in a waning block. Then the
    imaginary reflection index contributed by $i$ is
    \[
        (i,i+n)^*.
    \]
    In the reflection table, the diagonal entry is
    \[
        T[i,i] = [\,n\,]^*.
    \]
\end{itemize}

Following these rules, given the window notation of a TITO, its inversion set can be computed by analyzing every pair of residue classes. The resulting
inversion set is stored in a finite reflection table. The entry in row $a$ and
column $b$ records offsets $d$ such that the reflection index is
\[
(a,a+d),
\]
with $a+d \equiv b \pmod n$. If the entry is marked with a star, then it
represents the corresponding infinite family of reflection indices as defined above.


\paragraph{Example.}
The following example computes the inversion set of a TITO with period $n=4$.
\begin{verbatim}
from TITO_Explore.inversion_set import (
    tito_to_inversion_set,
    print_inversion_matrix,
    print_inversion_set,
)
from TITO_Explore.types import TranslationInvariantTotalOrder
tito = TranslationInvariantTotalOrder(
    n=4,
    num_blocks=2,
    vectors=[[0, 5, 6], [3]],
    waxing_waning=[0, 1],
)
print("TITO:", tito)
inversion_matrix = tito_to_inversion_set(tito)
print("\nInversion matrix:")
print_inversion_matrix(inversion_matrix)
print("\nInversion set:")
print_inversion_set(inversion_matrix)
\end{verbatim}
The output is:
\begin{verbatim}
TITO: TranslationInvariantTotalOrder(n=4, num_blocks=2, 
vectors=[[0, 5, 6], [3]], waxing_waning=[0, 1])
Inversion matrix:
         .        [1]        [2]          .
         .          .          .          .
         .          .          .          .
      [1]*       [2]*       [3]*       [4]*

Inversion set:
{ (0,1), (0,2), (3,4)*, (3,5)*, (3,6)*, (3,7)* }
\end{verbatim}

\subsection{Join Computation}
To compute the join (transitive closure) of two given TITOs with known inversion sets, we construct a unified directed graph. This graph encodes the combined inversion sets of both TITOs, allowing us to compute the complete inversion set of their join via path detection and edge weight updates.

As a preliminary step, we apply the weak-order comparison algorithm to the two given TITOs. If they are found to be comparable, the algorithm returns the larger TITO and terminates. Otherwise, it proceeds to the subsequent steps.

\paragraph{Setup: Graph representation and data structure.} 
Let $n$ be the period of the input TITOs. We define a directed graph with $n$ nodes, corresponding to the $n$ residue classes modulo $n$. Directed edges between nodes are assigned weights based on the inversion sets of the two TITOs to be joined. 
Computationally, this graph is represented by an $n \times n$ matrix, where each entry is a tuple $(B_{i,j}, c_{i,j})$ containing a basis set of smallest increments and an indicator for infinite growth, initially set to $0$. This initialization structure is illustrated in Table~\ref{tab:join_initialization}.

\begin{table}[htpb]
\centering
\caption{Initial matrix representation of inversion sets for TITO join computation. Each entry is a tuple $(B_{i,j}, c_{i,j})$ containing a list of the basis increments and an indicator for infinite growth, initialized to $0$.} 
\label{tab:join_initialization}
\renewcommand{\arraystretch}{1.5}
\begin{tabular}{@{}c|cccc@{}}
\toprule
$i \setminus j$ & $0$ & $1$ & $\cdots$ & $n-1$ \\ \midrule
$0$             & $(B_{0,0}, 0)$   & $(B_{0,1}, 0)$   & $\cdots$ & $(B_{0,n-1}, 0)$   \\
$1$             & $(B_{1,0}, 0)$   & $(B_{1,1}, 0)$   & $\cdots$ & $(B_{1,n-1}, 0)$   \\
$\vdots$        & $\vdots$         & $\vdots$         & $\ddots$ & $\vdots$           \\
$n-1$           & $(B_{n-1,0}, 0)$ & $(B_{n-1,1}, 0)$ & $\cdots$ & $(B_{n-1,n-1}, 0)$ \\ \bottomrule
\end{tabular}
\end{table}

\paragraph{Step 1: Set up initial matrix.} 
For every finite reflection index $\boldsymbol{(}a,b\boldsymbol{)}$ in the initial inversion sets, where $a \equiv i \pmod n$ and $b \equiv j \pmod n$, we append the increment $(b-a)$ to the basis increment list $B_{i,j}$. If the inversion is instead represented by an infinite reflection index $\boldsymbol{(}a,b\boldsymbol{)}^* = \{\boldsymbol{(}a,\, b + kn\boldsymbol{)} \mid k \in \mathbb{Z}_{\geq0}\}$, we similarly append $(b-a)$ to $B_{i,j}$ and update the infinite growth indicator $c_{i,j} \gets 1$. Applying this rule to all reflection indices across both TITOs sets up the matrix.

\paragraph{Step 2: Path detection for basis increment updates.} 
Given the initialized matrix, the algorithm modifies the classical Floyd--Warshall algorithm \cite{floydwarshallalgorithm, warshall1962theorem} to obtain the transitive closure of basis increments for each edge. Specifically, for every residue class from $0$ to $n-1$ acting as an intermediate node $k$, the algorithm evaluates all pairs of start and end nodes $(i, j)$. The basis sets are iteratively updated via Cartesian addition:
$$B_{i,j} \gets B_{i,j} \cup \{ x + y \mid x \in B_{i,k},\ y \in B_{k,j} \}$$
To prevent trivial infinite compounding from self-loops during this phase, we strictly require $k \notin \{i, j\}$. Because the algorithm iterates through all triplets of residue classes $(k, i, j)$ exactly once while performing Cartesian additions on basis sets of maximum cardinality $W$, the overall time complexity for this step is bounded by $O(n^3 W^2)$.

\paragraph{Step 3: Cycle detection for indicator updates.} 
Once the basis increments along each edge have been obtained, the algorithm identifies paths capable of generating infinite inversions and updates their corresponding indicators. We first examine the diagonal of the matrix for self-loops. If a diagonal entry $B_{i,i}$ is non-empty, traversing this loop generates an infinite sequence of increments $w + kn$ for $k \in \mathbb{Z}_{>0}$. This corresponds directly to an infinite family of imaginary reflection indices. We update $c_{i,i} \gets 1$. Subsequently, for any valid edge $i \to j$ where $B_{i,j} \neq \emptyset$ and $i \neq j$, if either endpoint possesses an infinite self-loop ($c_{i,i} = 1$ or $c_{j,j} = 1$), the inversions along the path $i \to j$ can also compound infinitely. Consequently, we update $c_{i,j} \gets 1$.

Taken together, Algorithm~\ref{alg:find_j_inv_set} computes the complete inversion set of the join, given two TITOs.

\begin{algorithm}
\caption{Computing the Reflection Table of the Join of two TITOs}
\label{alg:find_j_inv_set}
\begin{algorithmic}[1]
\Require The combined finitely represented inversion set $I = I_1 \cup I_2$ of two TITOs of period $n$

\State Initialize $M$ as an $n \times n$ matrix where every entry is $(\emptyset, 0)$

\Statex \Comment{\textbf{Step 1: Set up initial matrix}}
\For{each reflection index $R \in I$}
    \If{$R$ is of the finite form $\boldsymbol{(}a, b\boldsymbol{)}$}
        \State $i \gets a \pmod n$
        \State $j \gets b \pmod n$
        \State Append $(b - a)$ to $B_{i,j}$
    \ElsIf{$R$ is of the infinite form $\boldsymbol{(}a, b\boldsymbol{)}^*$}
        \State $i \gets a \pmod n$
        \State $j \gets b \pmod n$
        \State Append $(b - a)$ to $B_{i,j}$
        \State $c_{i,j} \gets 1$
    \EndIf
\EndFor

\Statex \Comment{\textbf{Step 2: Path detection for basis increment updates}}
\For{$k=0$ to $n-1$}
    \For{$i=0$ to $n-1$}
        \For{$j=0$ to $n-1$}
            \If{$k \neq i$ \textbf{and} $k \neq j$}
                \State $B_{i,j} \gets B_{i,j} \cup \{ x + y \mid x \in B_{i,k},\ y \in B_{k,j} \}$
            \EndIf
        \EndFor
    \EndFor
\EndFor

\Statex \Comment{\textbf{Step 3: Cycle detection for indicator updates}}
\For{$i=0$ to $n-1$}
    \If{$B_{i,i} \neq \emptyset$}
        \State $c_{i,i} \gets 1$ \Comment{Identify infinite self-loops}
    \EndIf
\EndFor

\For{$i=0$ to $n-1$}
    \For{$j=0$ to $n-1$}
        \If{$B_{i,j} \neq \emptyset$ \textbf{and} ($c_{i,i} = 1$ \textbf{or} $c_{j,j} = 1$)}
            \State $c_{i,j} \gets 1$ \Comment{Propagate infinite property to paths}
        \EndIf
    \EndFor
\EndFor

\State \Return $M$
\end{algorithmic}
\end{algorithm}

\paragraph{Step 4: Convert from the reflection table to a TITO object.}
To convert the reflection table back into a TITO object, we first examine each pair of residue classes and determine whether they belong to the same waxing/waning block or to different blocks. Then, our algorithm determines the order of the blocks and the order of the numbers within each block. 

For each pair of distinct residue classes $\{i,j\}$, we inspect the reflection-table entries involving $i$ and $j$. Each such entry has a direction, such as from $i$ to $j$, and an indicator. We call an entry \emph{starred} if its indicator is $1$, corresponding to an infinite family of reflections.

The pair $\{i,j\}$ is classified according to the following rules:

\begin{itemize}
    \item \textbf{Different blocks.}
    If there is exactly one starred directed reflection between $i$ and $j$, then $i$ and $j$ belong to different blocks.

    \item \textbf{Same waxing block.}
    If the reflections between $i$ and $j$ are unstarred in both directions, then $i$ and $j$ belong to the same waxing block.

    \item \textbf{Same waning block.}
    If in both directions, the reflections are starred, then $i$ and $j$ belong to the same waning block.
\end{itemize}

Then, the algorithm loops through each pair of residue classes that are identified as belonging to different blocks to determine the order of the blocks using the following rules. For residue classes $i<j$, where $i$ and $j$ are in different blocks, if the directed reflection from $i$ to $j$ has a star, then the block containing $j$ comes before the block containing $i$. Otherwise, the block containing $i$ comes before the block containing $j$. This follows from Cases 1 and 2 discussed in \ref{findinversionset}. 

After determining 1) the residue classes in each block, 2) whether the block is waxing or waning, and 3) the order among blocks, it is left to determine the exact order of numbers within each block. Here, we discuss waxing and waning blocks. The main idea is to first recover the exact number from each residue class in the block, and then determine the correct order of these recovered numbers. 

Here we abuse the notation by treating $i$ and $j$ both as residue classes and numbers between $0$ and $n-1$.

\paragraph{Waxing Block.}
Within a waxing block, the number of reflection indices is finite. Take the
smallest residue class $i$ as an anchor. For every other residue class $j$ in
the same block, the algorithm first determines which representative of the
residue class $j$ appears in the normalized window. Denote this representative
by $t_j$, where
\[
    t_j \equiv j \pmod n.
\]
If there is no finite reflection index from $i$ to the residue class $j$, then
the algorithm sets
\[
    t_j=j.
\]
Otherwise, we follow either Case 3 or Case 4 as discussed in \ref{findinversionset}. 

Following Case 3, if the directed reflection entry from i to j in the reflection table is non-empty, let $m_j$ be the largest offset value in that entry. The algorithm recovers
\[
    t_j=i+m_j+n, 
\] and the local window notation is $[i,t_j]$, where $i<t_j$.

Following Case 4, if the directed reflection entry from $j$ to $i$ is non-empty, let $n_j$ be the smallest offset in that entry. The algorithm recovers 
\[
    t_j = n_j-n, 
\] and the local window notation is $[i,t_j]$, where $i>t_j$. 

After these representatives are determined, the algorithm reconstructs the
order inside the waxing block. It begins with the anchor $i$ as the first entry of the block, following the convention in the normalization algorithm. Then it inserts the remaining representatives one at a time. A representative $t_j$ is inserted into the unique position for which its order relative to the representatives already placed agrees with the finite reflection indices in the reflection table. In a valid reflection table, this insertion process produces a unique order.

\paragraph{Waning Block.}
Within a waning block, the number of reflection indices is cofinite. Here, we follow either Case 5 or Case 6, both with starred reflection indices. Similarly, we take the smallest residue class $i$ as an anchor and place it as the first element in the reconstructed block. 

If there is no finite reflection index from $i$ to the residue class $j$, then
the algorithm sets
\[
    t_j=j.
\]
Otherwise, we follow either Case 5 or Case 6 from \ref{findinversionset}.

Following Case 5, let the offset in the directed reflection entry from residue class $i$ to $j$ be $m_i$, then we recover \[
    t_j=m_i-n,
\] and the local window is $\underline{[i,t_j]}$, where $i<t_j$. 

Following Case 6, let the offset in the directed reflection entry from residue class $j$ to $i$ be $n_j$, then the algorithm recovers 
\[
t_j=n_j, 
\] and the local window notation is $\underline{[i,t_j]}$, where $i>t_j$. 

Then, we follow a similar idea to insert the recovered numbers into the correct position. Here, we utilize Case 5 and Case 6 from \ref{findinversionset}. Suppose we are inserting $a$ from residue class $i$, and we are comparing it
with an existing non-anchor value $b$ from residue class $j$. We want to decide whether $a$ should appear before or after $b$ in the waning block. Without loss
of generality, suppose $i<j$.

First, we translate the pair so that the representative of residue class $i$
is equal to $i$. If $a$ is the representative of residue class $i$, then we
subtract $a-i$ from both entries. Thus the pair becomes
\[
    x=i,
    \qquad
    y=b-(a-i).
\]
Now the local two-element waning window has the form
\[
    \underline{[x,y]},
    \qquad x\equiv i \pmod n,\quad y\equiv j \pmod n.
\]

There are two possibilities. If
\[
    x<y,
\]
then the pair is in Case 5. In this case, the starred reflection entry from $i$ to $j$ has offset
\[
    y+n-i.
\]
Then the local order agrees with
\[
    \underline{[x,y]}.
\]
Otherwise, the local order must be reversed.

On the other hand, if
\[
    x>y,
\]
then the pair is in Case 6. In this case, the starred reflection entry from $i$ to $j$ has offset
\[
    j-i.
\]
Then the local order agrees with
\[
    \underline{[x,y]}.
\]
Otherwise, the local order must be reversed.

Using this pairwise comparison rule, the algorithm inserts the recovered values one at a time. The anchor remains fixed as the first entry of the block. 

\paragraph{Example.}
The following example computes the join of two TITOs and reconstructs the resulting TITO from the join reflection table.

\begin{verbatim}
from TITO_Explore.join import compute_join_reflection_table, compute_tito_join_tito
from TITO_Explore.types import TranslationInvariantTotalOrder

tito1 = TranslationInvariantTotalOrder(
    n=2,
    num_blocks=1,
    vectors=[[0, 3]],
    waxing_waning=[0],
)

tito2 = TranslationInvariantTotalOrder(
    n=2,
    num_blocks=1,
    vectors=[[2, 1]],
    waxing_waning=[0],
)

table_rows = compute_join_reflection_table(tito1, tito2)
joined_tito = compute_tito_join_tito(tito1, tito2)

print("Reflection table rows:")
for row in table_rows:
    print(row)

print("\nJoin TITO:")
print(joined_tito)
\end{verbatim}

The output is:

\begin{verbatim}
Reflection table rows:
(0, 2, 1)
(0, 1, 1)
(1, 2, 1)
(1, 3, 1)

Join TITO:
TranslationInvariantTotalOrder(n=2, num_blocks=1, vectors=[[0, -1]], waxing_waning=[1])
\end{verbatim}

Thus the join is the waning TITO
\[
    \underline{[0,-1]}.
\]

\section{Acknowledgments}
We would like to express our sincere gratitude to Chi Dinh for her thoughtful discussions with us and Allen Graham Hart for helpful feedback that led to an improvement in the package. We are deeply appreciative of the mentorship of Dr. Grant Barkley and Mia Smith, whose guidance, insight, and encouragement were invaluable to our work. We are grateful as well to Dr. Alejandro Bravo-Doddoli and The Lab of Geometry at Michigan for this wonderful research opportunity.

\clearpage
\printbibliography
\end{document}